\documentclass[12pt]{article}
\usepackage{graphicx}

\usepackage{amsmath,amsfonts,amsthm,amssymb}
\newcommand{\homo}{{\mathrm H}}
\def\co{\colon\thinspace}

\usepackage{indentfirst}

\usepackage{latexsym}

\newtheorem{Theorem}{Theorem}[section]
\newtheorem{Lemma}[Theorem]{Lemma}
\newtheorem{Corollary}[Theorem]{Corollary}

\theoremstyle{definition}

\newtheorem{Notation}[Theorem]{Notation}
\newtheorem{Example}[Theorem]{Example}

\begin{document}

\title{On standard forms of $1$--dominations between knots with same Gromov volumes}

\date{}

\maketitle

\author{\begin{center} MICHEL BOILEAU, boileau@picard.ups-tlse.fr\\
{\it Laboratoire \'Emile Picard\\Universit\'e Paul Sabatier,
TOULOUSE Cedex 4, France}\and

\ \\YI NI,  yni@math.columbia.edu\\
{\it Department of Mathematics, Columbia University\\ MC 4406,
2990 Broadway, New York, NY 10027}\and

\ \\SHICHENG WANG, wangsc@math.pku.edu.cn\\
{\it LMAM, Department of Mathematics, \\Peking University, Beijing
100871 China}\end{center}}

\section{Introduction and notations}

All knots are in the $3$--sphere $S^3$. For basic terminologies in
knot theory and in $3$--manifold theory, see \cite{R}, \cite{He}
and \cite{Ja}.

We recall the following relation on the set of knots in $S^3$: let
$k_1$ and $k_2$ be two knots, we say $k_1 \geq k_2$, or
equivalently say that $k_1$ 1--dominates $k_2$, if there is a
proper degree $1$ map $f\co E(k_1)\to E(k_2)$, where $E(k_i)$ is
the knot exterior of $k_i$. If $k_1\ge k_2$ but $k_1\ne k_2$, we
often write $k_1>k_2$, or equivalently say that the  1--domination
is non-trivial.

Following the classical results of \cite{Wal} and \cite{GL}, it is
known that the relation $\ge$  is a partial order on knots in
$S^3$.

In general, when $k_1\ge k_2$, the relation of $k_1$ and $k_2$ is
not known, and there is no fine description of the degree 1 map,
up to homotopy, realizing the 1--domination $k_1>k_2$. Recall that
a simplest and a most common construction of 1--domination of
$k_1>k_2$ is to choose $k_1$ to be a satellization of $k_2$ and
$f$ realizing the 1--domination to be the de-satellization, and on
the other hand there are many sophisticated constructions, see
\cite{Ka}, \cite{Ru}, \cite{BW1}, \cite{BW2}, \cite{BNW},
 \cite{ORS} and so on.

In this note we show  that 1--domination between knots is
de-satellization under certain conditions.

\begin{Theorem}\label{pinching}
Suppose that any companion of $k$ is prime. If $k\ge k'$ with the
same Gromov volume, then $k'$ can be obtained from $k$ by finitely
many de-satellizations.
\end{Theorem}

The condition of ``same Gromov volume" clearly can not be removed,
according to the constructions in the papers mentioned above. We
will also give a new construction of 1--domination between knots
with same Gromov volume to show that the condition ``any companion
of $k$ is prime" can not be removed.

The corollary below supports a general opinion that the
1-domination partial order reflects the complexity of knots (see a
survey \cite{W}). By a theorem of Schubert \cite{Sc}, we have

\begin{Corollary}
Suppose that any companion of $k$ is prime. If $k>k'$ with the
same Gromov volume, then $b(k)>b(k')$, where $b$ is the bridge
number.
\end{Corollary}

The paper is organized as follows. After listing some known useful
facts, a general study of maps between Seifert pieces and graph
pieces in knot complements is given in \S 2,
Theorem~\ref{pinching} will be proved in \S 3, and the new
construction of 1-domination will be given in \S 4. Below we will
fix some notions for the remaining sections.

\begin{Notation}
For each solid torus in $S^3$, we specify its longitude to be the
one which is homologous to zero in the complement. Let $k_1$ be a
geometrically essential knot \cite[p110]{R} in an unknotted solid
torus $V\subset S^3$ and $k_2$ be another knot. Let $h\co V\to
N(k_2)$ be a longitude preserving homeomorphism, then the new knot
$k=h(k_1)$ is called the satellite knot of $k_2$, and $k_2$ is a
companion of $k$.

The reversing process of satellization,   given  by pinching
$E(k_2)$, the exterior of the companion  to a solid torus,
produces a proper degree one map  $f\co E(k) \to E(k_1)$, which
will be called a {\it de-satellization}.
\end{Notation}

\begin{Notation} Let $T(p_1,q_1;p_2,q_2;\dots;p_n,q_n)$ be the
iterated torus knot, which is the $(p_1,q_1)$--cable of the
$(p_2,q_2)$--cable of \dots the $(p_n,q_n)$--torus knot. (When we
say ``$(p,q)$--cable", $p$ denotes the winding number.) The
exterior of the knot is denoted by
$E=E(p_1,q_1;p_2,q_2;\dots;p_n,q_n)$. Let
$C=C(p_1,q_1;p_2,q_2;\dots;p_n,q_n)$ denote the ``iterated cable
space", that is, $E$ with an open neighborhood of the singular
fiber corresponding to $(p_n,q_n)$ removed. $E$ is a graph
manifold, the Seifert pieces are denoted by $C(p_1,q_1)$,
\dots,$C(p_{n-1},q_{n-1})$,$E(p_n,q_n)$, $\partial E=T_0$, the JSJ
tori are denoted by $T_1,\dots,T_{n-1}$, where $\partial C(p_i,
q_i)=T_{i-1}\sqcup T_i$. $C$ is also a graph manifold, the Seifert
pieces are $C(p_1,q_1)$,\dots,$C(p_n,q_n)$, $\partial C=T_0\sqcup
T_n$, the JSJ tori are denoted by $T_1,\dots,T_{n-1}$, where
$\partial C(p_i, q_i)=T_{i-1}\sqcup T_i$. Suppose $\alpha$ is a
slope on $T_n$, then $C(\alpha)=C(p_1,q_1;\dots;p_n,q_n;\alpha)$
denotes the manifold obtained by Dehn filling along $\alpha$.

$E$ and $C$ are submanifolds of $S^3$, $T_i$ bounds a solid torus
$K_i$ in $S^3$. Suppose $\mu_i\subset T_i$ is the meridian of
$K_i$, and $\lambda_i\subset T_i$ is the longitude.
\end{Notation}

\begin{Notation}\label{NotaPm}
Let $D_0$ be a disc and $D_1, ... , D_n$ be sub-discs in the
interior of $D_0$, and denote $\partial D_i$ by $c_i$, and the
$n$-punctured disc $D_0\setminus \cup_{i=1}^n D_i$ by $P_n$. Then
$\partial P_n=\cup_{i=0}^n c_i$. Note that $P_1$ is an annulus.
Once $D_0$ is oriented, then $P_n$ and all $c_i$ are oriented.
\end{Notation}

\begin{Notation}
Let $f\co M\to N$ be a map between orientable compact connected
$n$--manifolds. We say that $f$ is {\it proper} if
$f^{-1}(\partial N)=\partial M$. We say that $f$ is {\it
allowable} if $f$ is proper and the degree of all possible
restrictions $f|\co F\to S$ have the same sign, where $F$ is a
component of $\partial M$ and $S$ is a component of $\partial N$.
\end{Notation}

\noindent{\bf Acknowledgements.} This research was conducted
during the periods the second author was a graduate student at
Princeton University and was employed by the Clay Mathematics
Institute as a Liftoff Fellow. The second author was partially
supported by a Graduate School Centennial Fellowship at Princeton
University. The third author is partially supported by grant
No.10631060 of the National Natural Science Foundation of China.

\section{Proper maps between Seifert pieces and graph pieces in knot complements}

The following four known facts, see \cite{Go}, \cite{Ja}, \cite{Ro}
and \cite{So} respectively, will be repeatly used in this paper.

\begin{Lemma}{\bf\cite{Go}}
In $C(p_i,q_i)$, we have the following relations in homology:
$$p_i[\mu_{i-1}]=[\mu_i],\quad[\lambda_{i-1}]=p_i[\lambda_i].$$
Moreover the regular Seifert fiber of $C(p_i,q_i)$ is homologous
to $p_iq_i[\mu_{i-1}]+[\lambda_{i-1}]$ on $T_{i-1}$, and
homologous to $q_i[\mu_i]+p_i[\lambda_i]$ on $T_i$.
\end{Lemma}

\begin{Lemma}{\bf\cite{Ja}}
Let $P$ be a Seifert piece of the JSJ decomposition of $E(k)$. Then
 $P$ is either $E(p,q)$, or $C(p,q)$, or $P_m\times S^1$.
Moreover $P_m\times S^1, m>1$ appears if and only some companion
of $k$ is not prime.
\end{Lemma}

\begin{Lemma}{\bf\cite{Ro}} Let $f\co M\to N$ be an allowable degree 1 map between
aphereical Seifert manifolds. Then $f$ is homotopic to a fiber
preserving pinch.
\end{Lemma}

\begin{Lemma}{\bf\cite{So}}
If $f\co M\to N$ is a proper map of degree $d$ between Haken
manifolds such that $||M||=d||N||$, then $f$ can be homotoped to
send $H(M)$ to $H(N)$ by a covering, where $||*||$ is the Gromov
norm and $H(*)$ is the hyperbolic part under the JSJ decomposition.
\end{Lemma}

Below  we give some general study of maps between Seifert pieces
and graph pieces in knot complements.

\begin{Lemma}\label{rong}
Any proper degree 1 map $f: E(p,q)\to E(p',q')$ between torus knot
complements is homotopic to a homeomorphism.
\end{Lemma}
\begin{proof}
The lemma is known since that all torus knots are minimal (see
\cite{BW1}). It is also a direct corollary of \cite{Ro}: Since each
manifold involved has only one boundary component, $f$ is an
allowable degree 1 map. Since each Seifert manifold involved has a
unique Seifert fibration, then by \cite{Ro}, $f$ is homotopic to a
fiber preserving pinch. Since any non-trivial pinch will decrease
either the genus of the orbifold, or the number of singular fibers,
and since both the genus of the orbifold and the number of singular
fibers of $E(p,q)$ and $E(p',q')$ are the same, the pinch must be
trivial, therefore the lemma is verified.
\end{proof}

\begin{Lemma}\label{fiberpre}
Suppose $M$ is a Seifert manifold with a $\pi_1$--injective
boundary component $T$  and $f\co C(p_1,q_1)\to M$ is a proper map
such that $f|\co T_0\to T$ is a homeomorphism.  Let
$t_1\in\pi_1(C(p_1,q_1))$ and $t\in \pi_1(M)$ represent regular
fibers of the corresponding Seifert manifolds. Then the following
statements hold.

(1) $f_*(\pi_1(C(p_1,q_1)))$ is not an abelian group;

(2) $f_*(t_1)=t^{\pm1}$ if $M$ has a unique Seifert fibration.
\end{Lemma}
\begin{proof}
Pick a base point of $C(p_1,q_1)$ in $T_0$, and a base point of
$M$ in $T$. Then $\pi_1(T_0)$ is naturally a subgroup of
$\pi_1(C(p_1,q_1))$, and $\pi_1(T)$ is naturally a subgroup of
$\pi_1(M)$.

Assume $f_*(\pi_1(C(p_1,q_1)))$ is an abelian group. Since $f|T_0$
is a homeomorphism, and $\pi_1(T)$ is a maximal abelian subgroup
of $\pi_1(M)$, $f_*(\pi_1(C(p_1,q_1)))$ must be $\pi_1(T)$.
Moreover,
$$f_*\co\pi_1(C(p_1,q_1))\to\pi_1(T)$$
factors through $\homo_1(C(p_1,q_1))$. $\lambda_0,\lambda_1$
represent elements in $\pi_1(C(p_1,q_1))$, $\lambda_0$ is the
$p_1$--multiple of $\lambda_1$ in $\homo_1(C(p_1,q_1))$, but
$f_*(\lambda_0)$ is a primitive element in $\pi_1(T)$, we get a
contradicition.

Since $t_1$ commutes with $\pi_1(C(p_1,q_1))$, and
$f_*(\pi_1(C(p_1,q_1)))$ is non-abelian, $f_*(t_1)$ must be a
power of $t$. Since $f\co T_0\to T$ is a homeomorphism,
$f_*(t_1)=t^{\pm1}$.
\end{proof}

\begin{Lemma}\label{onecomp}
Let
$$f\co C(\alpha)=C(p_1,q_1;\dots;p_n,q_n;\alpha)\to E(p,q)$$
be a proper map, and the restriction of $f$ to $T_0$ is a
homeomorphism. Then the restriction of $f$ to $T_1$ is not
$\pi_1$--injective.
\end{Lemma}
\begin{proof}
Pick a basepoint $b$ of $C(\alpha)$, $b\in T_0$, choose a simple
curve $\gamma$ connecting $b$ to $T_{n-1}$, such that $\gamma\cap
T_i$ consists of a single point. Let $\gamma\cap T_i$ be the base
point in $T_i$ and $E(p_{i+1},q_{i+1})$. Using a path on $\gamma$,
we can view $\pi_1(T_i)$ and $\pi_1(E(p_{i+1},q_{i+1}))$ as
subgroups of $\pi_1(C(\alpha))$. Let
$f_*\co\pi_1(C(\alpha))\to\pi_1(E(p,q))$ be the induced map on
$\pi_1$. Let $T_0'=\partial E(p,q)$.

Let $t_i\subset \pi_1(C(p_i,q_i))$ and $t\subset\pi_1(E(p,q))$
represent the regular Seifert fibers in the corresponding Seifert
manifolds. By Lemma \ref{fiberpre}, we can assume $f_*(t_1)=t$.

If $n=1$, then the conclusion trivially holds (since $\alpha$ is
in the kernal), so we assume $n>1$. The element $t_1$ is contained
in $\pi_1(T_1)$. In fact, $t_1$ is homologous to
$q_1[\mu_1]+p_1[\lambda_1]$ in $T_1$. Let $x$ denote $f_*(\mu_1)$.
Assume the restriction of $f$ on $T_1$ is $\pi_1$--injective, then
$x,t$ generate a $\mathbb Z\oplus\mathbb Z$--subgroup of
$\pi_1(E(p,q))$.

The fiber $t_2$ is homologous to $p_2q_2[\mu_1]+[\lambda_1]$ on
$T_1$, hence not a power of $t_1$ in $\pi_1(T_1)$. So $f_*(t_2)$
is not a power of $t$. But $t_2$ commutes with
$\pi_1(C(p_2,q_2))$, so $f_*(\pi_1(C(p_2,q_2)))$ is an abelian
group. Hence
$$f_*\co\pi_1(C(p_2,q_2))\to\pi_1(E(p,q))$$
factors through $\homo_1(C(p_2,q_2))$.

In $C(p_2,q_2)$, $p_2(q_1[\mu_1]+p_1[\lambda_1])$ is homologous to
$q_1[\mu_2]+p_1p_2^2[\lambda_2]$, hence the corresponding element
in $\pi_1(T_2)$ is mapped by $f_*$ to $t^{p_2}$. By the same
reason, $f_*(\mu_2)=x^{p_2}$. So $f|T_2$ is $\pi_1$--injective.
Moreover, $t_3$ is homologous to $p_3q_3[\mu_2]+[\lambda_2]$ in
$T_2$, it is linearly independent with
$q_1[\mu_2]+p_1p_2^2[\lambda_2]$, since $\gcd(p_1,q_1)=1$. Hence
$f_*(t_3)$ is not a power of $t$. But $t_3$ commutes with
$\pi_1(C(p_3,q_3))$, so $f_*(\pi_1(C(p_3,q_3)))$ is an abelian
group.

Argue as before, we find that $f_*(\mu_3)=x^{p_2p_3}$, and the
loop corresponding to $q_1[\mu_3]+p_1p_2^2p_3^2[\lambda_3]$ on
$T_3$ is mapped to $t^{p_2p_3}$ by $f_*$. Hence $f|T_3$ is
$\pi_1$--injective, and $f_*(t_4)$ is not a power of $t$.

Go on with such argument, we finally show that $f|T_{n-1}$ is
$\pi_1$--injective, and $f_*(t_n)$ is not a power of $t$, where
$t_n$ represents the regular fiber of $C(p_n,q_n)$. Thus
$f_*(\pi_1(C(p_n,q_n)))$ is an abelian group, and therefore the
group $f_*(\pi_1(C(p_n,q_n;\alpha)))$ is also abelian. Then
$f_*|\pi_1(C(p_n,q_n;\alpha))$ factors through
$\homo_1(C(p_n,q_n,\alpha))\cong \mathbb Z\oplus \mathbb Z_b$ for
some positive integer $b$, which contradicts to the fact that
$f|T_{n-1}$ is $\pi_1$--injective.
\end{proof}

\begin{Lemma}\label{twocomp}
Suppose $\partial C(p,q)=T_0'\sqcup T_1'$, where $T_0'$ bounds a
neighborhood of the torus knot $T(p,q)$, and
$$f\co C(p_1,q_1;\dots;p_n,q_n)\to C(p,q)$$
is a proper map.

(1) If $n>1$, then $f$ cannot map $T_0$ homeomorphically to $T_0'$.

(2) If $n=1$, and $f$ maps $T_0$ homeomorphically to $T_0'$, then
$f$ is homotopic to a homeomorphism.
\end{Lemma}
\begin{proof}
Assume $f$ maps $T_0$ homeomorphically to $T_0'$. We claim that
$f(T_n)=T_1'$. Otherwise $f(T_n)=T_0'$. Let $f_{\#}$ be the
induced map on homology. $f_{\#}([\lambda_n])$ is an integral
linear combination of $f_{\#}([\mu_0])$ and $f_{\#}([\lambda_0])$,
but $[\lambda_n]$ is equal to $\frac1P[\lambda_0]$, where
$P=p_1p_2\dots p_n$. We get a contradiction.

Now $f(T_n)=T_1'$. Since $f|T_0$ is a homeomorphism, $\text{deg}\:
f=\text{deg}\: f|T_n=1$. We can homotope $f$, so that $f|T_n$ is a
homeomorphism. Moreover,
$$f_{\#}\co\homo_1(C)\to\homo_1(C(p,q))$$
is an isomorphism.

By Lemma \ref{fiberpre}, we can assume $f_*(t_1)=f_*(t_n)=t$. In
$\homo_1(C)$, we have
$$[t_1]=p_1q_1[\mu_0]+[\lambda_0]=p_1q_1[\mu_0]+P[\lambda_n],$$
$$[t_n]=q_n[\mu_n]+p_n[\lambda_n]=q_nP[\mu_0]+p_n[\lambda_n].$$
Since $f_{\#}$ is an isomorphism and $[\mu_0],[\lambda_n]$
generate $\homo_1(C)$, we must have
$$p_1q_1=q_n P, \,\,\, P=p_n.$$

If $n>1$, it is impossible since  $p_1>1$.

If $n=1$, then we have a proper allowable degree map $f\co
C(p_1,q_1)\to C(p,q)$. Applying Rong's result as in the proof of
Lemma~\ref{rong}, one shows that $f$ is homotopic to a
homeomorphism.
\end{proof}

\begin{Lemma}\label{pants} Let $M$ be either
$E(p_1,q_1;\dots;p_n,q_n)$ or $C(p_1,q_1;\dots;p_n,q_n)$, let
$P_m$ denote the $m$-punctured disk, where $m>0$. Then there is no
proper map $f\co M\to P_m\times S^1$ such that $f$ restricts to a
component of $\partial M$ is a homeomorphism.
\end{Lemma}
\begin{proof} Assume $f$ maps $T_0$ homeomorphically to $T_0'$,
a component of $\partial P_m\times S^1$.

If $M$ is a knot space, then $[T'_0]=f_\#([T_0])$ is null
homologous in $P_m\times S^1$, which implies $m=0$, a
contradiction.

Now suppose that $M$ is an iterated cable space with boundary
$T_0$ and $T_n$. Since $[\lambda_{0}]=p[\lambda_n]$ in
$\homo_1(M;\mathbb Z)$, where $p=p_1... p_n>1$, we have
$f_{\#}([\lambda_{0}])=p f_{\#}([\lambda_n])$ in
$\homo_1(P_m\times S^1;\mathbb Z)=\mathbb Z^{m+1}$. There are two
subcases:

(a) $f_{\#}([T_n])=k[T_0']$, $k\in\mathbb Z$;

(b) $f_{\#}([T_n])=k[T_1']$, $k\in\mathbb Z$, $T_1'\ne T_0'$,
$T_1'$ is a component of $\partial P_m\times S^1$.

In the subcase (a), since $[T_0]+[T_n]=0$, we have
$(k+1)[T_0']=0$, which implies that $k=-1$.  Now both
$f_{\#}([\lambda_{0}])$ and $f_{\#}([\lambda_{n}])$ are homologous
to
 closed curves on $T_0'$, and in particular
$f_{\#}([\lambda_{0}])$ is a primitive element in
$\homo_1(T'_0;\mathbb Z)=\mathbb Z^{2}$. Note that
$f_{\#}([\lambda_{0}])=p f_{\#}([\lambda_n])$ in
$\homo_1(P_m\times S^1;\mathbb Z)$, and the homomorphism
$\homo_1(T'_0;\mathbb Z)\to \homo_1(P_m\times S^1;\mathbb Z)$
induced by the inclusion is injective, so $f_{\#}([\lambda_{0}])=p
f_{\#}([\lambda_n])$ in $\homo_1(T'_0;\mathbb Z)$, which is
impossible since $f_{\#}([\lambda_0])$ is primitive.

In the subcase (b), since $[T_0]+[T_n]=0$, we have
$[T'_0]+[T'_1]=0$, which is impossible if $m>1$. If $m=1$, then
$P_1\times S^1= T_0'\times [0,1]$, and the homomorphism
$\homo_1(T'_0;\mathbb Z)\to \homo_1(T_0'\times [0,1];\mathbb Z)$
induced by the inclusion is an isomorphism, and again
$f_{\#}([\lambda_{0}])=p f_{\#}([\lambda_n])$ in
$\homo_1(P_1\times S^1;\mathbb Z)$, which is impossible.

In either case we reach a contradiction.
\end{proof}

\section{Proof of Theorem~\ref{pinching}}

The dual graph $\Gamma(k)$ to the JSJ-decomposition of $E(k)$ is a
rooted tree, where the root is corresponding to the unique vertex
manifold containing  $\partial E(k)$. Let $\Gamma_0(k) \subset
\Gamma(k)$ be the maximal connected subtree which contains the
root such that the restriction of $f$, up to homotopy,  to the
connected submanifold $M(\Gamma_0)$ associated to $\Gamma_0$ is a
homeomorphism to its image, and moreover the restriction of $f$ to
each leaf torus of $\Gamma_0$ is $\pi_1$-injective.

Since $k$ and $k'$ have the same Gromov volume, by \cite{So}, $f$
can be homotoped so that $f$ maps the hyperbolic pieces of $E(k)$
homeomorphically to the hyperbolic pieces of $E(k')$.

If $f\co E(k)\to E(k')$ is homotopic to a homeomorphism, then
Theorem~\ref{pinching} is automatically true. So below we assume
that $f$ is not homotopic to a homeomorphism. Then $M(\Gamma_0)\ne
E(k)$.

Let $T_0$ be the torus corresponding to a leaf of $\Gamma_0$, and
$X_0$ ($\not\subset M(\Gamma_0)$) be the JSJ piece adjacent to
$T_0$. Then $X_0$ must be a Seifert piece. Since $f|T_0$ is
$\pi_1$--injective, $f|X_0$ is non-degenerate, and it follows that
we can push $f(X_0)$ into a Seifert piece $X'_0$ of the JSJ
decomposition of $E(k')$. Let $T'_0=f(T_0)\subset\partial X'_0$,
then $f|: T_0\to T'_0$ is a homeomorphism. By the definition of
$\Gamma_0(k)$, we have a JSJ piece $X\ne X_0$ of $E(k)$ adjacent
to $T_0$ such that $f|\co X\to X'$ is a homeomorphism, where $X'$
is a JSJ piece of $E(k')$ adjacent to $T'_0$.

Let $U$ be the maximal connected graph submanifold of $E(k)$  such
that $X_0\subset U$ and $T_0\subset
\partial U$. Since we assume that any companion of $k$ is prime, then $U$ is in the
form of either  $E(p_1,q_1;\dots;p_n,q_n)$ or
$C(p_1,q_1;\dots;p_n,q_n)$.

\begin{Lemma}\label{non-torus-knot} $U\ne E(p_1,q_1)$, hence $T_1\ne \emptyset$.
\end{Lemma}
\begin{proof} Otherwise we have  $U=X_0= E(p_1,q_1)$ and then
$f(T_0)=T'_0$ is homologous to zero in $X_0'$, which implies
$\partial X'_0=T'_0$ and therefore $X'_0=E(p',q')$. Then we have
map $f|: E(p,q)\to E(p',q')$ which is degree 1 on the boundary,
and therefore degree 1 itself. By Lemma \ref{rong}, $f|$ is
homotopic to  a homeomorphism, and therefore contradicts to the
maximality of $\Gamma_0$.
\end{proof}

Below we name JSJ-tori in $U$ after $T_1$ as $T_2,..., T_n$ in
order.

\begin{Lemma}\label{non-injective1} $f|T_i$ is not $\pi_1$--injective for some $i$.
\end{Lemma}
\begin{proof}Otherwise the restriction of $f$ to any Seifert
piece in $U$ is non-degenerate. By homotoping $f$, we can assume
$f^{-1}(X_0')$ is the union of some Seifert pieces in $E(k)$.

Let $G$ be a component of $f^{-1}(X_0')$ containing $X_0$. The $G$
is either $E(p_1,q_1;\dots;p_l,q_l)$ or
$C(p_1,q_1;\dots;p_l,q_l)$.

\noindent{\bf Claim 1.\;} $X'_0=E(p',q')$, and $X'\ne X'_0$.

\begin{proof} By Lemma \ref{pants}, $X_0'$ is not $P_m\times S^1$, $m\ge 1$. Hence either
$X'_0=C(p',q')$ or $X'_0=E(p',q')$.

Suppose first $X'_0=C(p',q')$. By simple homological reason $G$
cannot be $E(p_1,q_1;\dots;p_l,q_l)$. By Lemma \ref{twocomp}, $G$
cannot be $C(p_1,q_1;\dots;p_l,q_l)$, $l>1$; moreover if
$C=C(p_1,q_1)$, then $f|\co C(p_1, q_1)\to C(p',q')$ is homotopic
to a homeomorphism, which contradicts to the maximality of
$\Gamma_0$.

Hence  $X'_0=E(p',q')$. Since $X'$, which is homeomorphic to $X$,
has at least two boundary components, we have $X_0'\ne X'$.
\end{proof}

\noindent{\bf Claim 2.\;} $f^{-1}(X_0')\cap U=U$.

\begin{proof}
Let $S'$ be a Seifert surface of $E(p',q')$. Since $f|\co T_0\to
T'_0$ is a homeomorphism, up to a homotopy relative to $T_0$,  we
may assume that $f^{-1}(S')$ is incompressible, and moreover there
is only one component of $f^{-1}(S')$, denoted by $S$, with
$\partial S$ a circle $c$. Since $f(X)= X'$, $X'\ne X'_0$, it
follows $f^{-1}(S')\cap int X=\emptyset$. Since $T_0$ is
separating and $S$ is connected, we must have $S\subset E(k_0)$,
hence $c=\lambda_0$, where $E(k_0)$ is a component separated by
$T_0$ containing $U$. Since the winding number of each JSJ torus
$T_i$ is non-zero with respect to $T_0$, we have $S\cap T_i\ne
\emptyset$ for each $i$, and it follows that $f^{-1}(X_0')\cap
U=U$.
\end{proof}

\noindent{\bf Claim 3.\;} $U=E(p_1,q_1;\dots;p_n,q_n)$.

\begin{proof}
If $U=C(p_1,q_1;\dots;p_n,q_n)$. Let $Y$ be the JSJ piece of
$E(k)-U$ adjacent to $T_n$. By the definition of $U$, $Y$ must be
a hyperbolic piece, so $f|Y$ must be a homeomorphism. Since
$f(T_n)\subset T'_0$, we must  have $f(Y)\subset X'$ and it
implies that $X'$ is a hyperbolic piece. Since $f\co X\to X'$ is
homeomorphism by our assumption, it follows that $X$ is a
hyperbolic piece. Therefore $f$ send two different hyperbolic JSJ
pieces of $E(k)$ to one hyperbolic JSJ piece of $E(k')$, it
contradicts that $f|$ on the hyperbolic part is a homeomorphism.
\end{proof}

Now we have a proper map $f\co E(p_1,q_1;\dots;p_n,q_n)\to
E(p',q')$ which is a homeomorphism on the boundary. By Lemma
\ref{onecomp}, $f|T_1$ is not $\pi_1$--injective, which
contradicts to the assumption we made before.

This finishes the proof of Lemma \ref{non-injective1}.
\end{proof}

\begin{Lemma}\label{non-injective2} $f|T_1$ is not $\pi_1$--injective.
\end{Lemma}
\begin{proof}By Lemma~\ref{non-injective1},
some $f|T_i$ is not $\pi_1$--injective for $T_i$ in $U$. We may
assume that $f|$ is $\pi_1$--injective on $T_i$ for $i<k$ and that
$f|$ is not $\pi_1$--injective on $T_k$. We have
$f(C(p_1,q_1;...;p_k,q_k))\subset X'_0$. Since $f|T_k$ is not
$\pi_1$--injective, there is a simple loop $\alpha\in T_i$ in the
kernel of $f_*$. Therefore we get a map $f|\co
C(p_1,q_1;...;p_k,q_k;\alpha)\to X'_0$ such that $f|T_0$ is a
homeomorphism. A homological argument shows that $X'_0=E(p,q)$. By
Lemma~\ref{onecomp} $f|T_1$ is not $\pi_1$--injective.
\end{proof}

\begin{proof}[Proof of Theorem~\ref{pinching}] Let $V=M(\Gamma_0)$,
$V'=f(V)$. Then $f|\co V\to V'$ is a homeomorphism. Denote the
knot complement separated by $T_i$ in $E(k)$ by $E(k_i)$, $i=0,1$
and $W=E(k)\setminus E(k_0)$. Then we have
$E(k_0)=C(p_1,q_1)\cup_{T_1}E(k_1)$ and there is a proper degree
one map
$$f\co E(k)= W \cup_{T_0}C(p_1,q_1)\cup_{T_1}E(k_1)\to E(k')$$
such that $f(C(p_1,q_1))\subset X'_0$, $f|\co T_0\to T'_0$ is a
homeomorphism, and a simple closed curve $\alpha\subset T_1$ lies
in the kernel of $f|T_1$. Then the proper degree one map $f\co
E(k)\to E(k')$ induces a factorization
$$E(k)\longrightarrow W
\cup_{T_0}C(p_1,q_1;\alpha)\cup_{\alpha^*}E(k_1,
\alpha)\stackrel{\hat f}{\longrightarrow} E(k') .\leqno(1)$$  Here
$C(p_1,q_1;\alpha)$ and $E(k_1, \alpha)$ are $3$--manifolds
obtained by Dehn filling along $\alpha\subset T_1$ on $C(p_1,q_1)$
and $E(k_1)$ respectively and
$C(p_1,q_1;\alpha)\cup_{\alpha^*}E(k_1, \alpha)$ is obtained by
identifying the core $\alpha^*$ of filling solid tori in
$C(p_1,q_1;\alpha)$ and $E(k_1,\alpha)$.

Since $E(k_1,\alpha)$ is a closed $3$--manifold, it makes no
contribution to the degree of the proper degree one map $f$ and we
have
$$\hat f|\co  W \cup_{T_0}C(p_1,q_1;\alpha)\to E(k')\leqno(2)$$
is a proper degree one map. Since
$$||E(k')||=||E(k)||\ge || W \cup_{T_0}C(p_1,q_1)||
\ge ||W \cup_{T_0}C(p_1,q_1;\alpha)||\ge ||E(k')||,$$ we have
$||E(k)||= || W \cup_{T_0}C(p_1,q_1)||$ and therefore $E(k_1)$ is a
graph manifold, it follows that
$$C(p_1,q_1)\cup_{T_1}E(k_1)=C(p_1,q_1)\cup_{T_1}E(p_2,q_2;...;p_n,q_n)\leqno(3)$$
where
$C(p_1,q_1)\cup_{T_1}E(p_2,q_2;...;p_n,q_n)=E(p_1,q_1;...;p_n,q_n)$

Moreover since $\hat f(C(p_1,q_1; \alpha))\subset X'_0$, $f|\co
T_0\to T'_0$ is a homeomorphism, it follows that $X_0'=E(p',q')$ and
$$\hat f|\co
C(p_1,q_1; \alpha) \to E(p',q')\leqno(4)$$ is homotopic to a
homeomorphism. Finally we have
$$f\co  W \cup_{T_0}C(p_1,q_1)\cup_{T_1}E(p_2,q_2;...;p_n,q_n)\to  E(k')=W'
\cup_{T'_0}E(p',q').\leqno(5)$$

Let $S'$ be a Seifert surface of $E(p',q')$, then up to a homotopy
relative to $T_0$,  we may assume that $f^{-1}(S')$ is
incompressible, and moreover there is only one component of
$f^{-1}(S')$, denoted by $S$, with $\partial S$ a circle. Let $X$
be a JSJ piece of $E(k)$ adjacent to $X_0$ along $T_0$, and let
$X'$ be a JSJ piece of $E(k')$ adjacent to $X'_0$ along $T'_0$. By
our choice of $T_0$, $f|X$ is a homeomorphism. Since $X$ has at
least two boundary components while $X'_0$ has only one boundary
component, we must have $f(X)\subset X'$ and therefore
$f^{-1}(S')\cap int X=\emptyset$. Since $T_0$ is separating and
$S$ is connected, we must have $S\subset
E(p_1,q_1;p_2,q_2;...;p_n,q_n)$ and therefore it is a Seifert
surface of $E(p_1,q_1;p_2,q_2;...;p_n,q_n)$ which intersects $T_1$
in parallel copies of $\lambda_1$. It follows that
$\alpha=\lambda_1$. Now we rewrite (1) as
$$E(k)\longrightarrow V \cup_{T_0}C(p_1,q_1;\lambda_1)\cup_{\lambda_1^*}E(k_1, \lambda_1)
\stackrel{\hat f}{\longrightarrow}
W'\cup_{T'_0}E(p',q').\leqno(6)$$

Note that the core $\lambda^*_1$ of the filling solid torus is a
retractor of $E(k_1, \lambda_1)$, and $\hat f|\co
C(p_1,q_1;\lambda_1)\to E(p',q')$ is homotopic to a homeomorphism
by \cite{Ro}. Now we have a further factorization
\begin{eqnarray*}
E(k)&\to& W \cup_{T_0}C(p_1,q_1;\lambda_1)\cup_{\lambda_1^*}E(k_1,
\lambda_1)\\ &\to& W \cup_{T_0}C(p_1,q_1;\lambda_1)=W
\cup_{T_0}E(p',q')\to W' \cup_{T'_0}E(p',q').
\end{eqnarray*}

Hence $f$ factors through the de-satellization:
$$E(k)\to W \cup_{T_0}E(p',q')\to E(k').$$

Clearly $W\cup_{T_0}E(p',q')=E(k'')$ for some knot $k''$ in $S^3$.
Moreover any companion of $k''$ is prime, $k''$ and $k'$ have the
same simplicial volume. So we can repeat the above process to
degree one map $E(k'')\to E(k')$. Since any knot admits at most
finitely many de-satellization, we finish the proof of
Theorem~\ref{pinching}.\end{proof}

\section{New construction}

\begin{Example} We construct a degree one map from a graph knot (i.e., the complement of the
knot is a graph manifolds) to a torus knot which is not a
de-satellization.

\vskip 0.5 truecm

Below  $c_i$ and $P_n$ are given in Notation~\ref{NotaPm}. We use
$\bar T(3,2)$ to denote the mirror image of $T(3,2)$ and $\bar
E(3,2)$ to donote the exterior of $\bar T(3,2)$.

\begin{Lemma}[Schubert]\label{Schubert}
 The JSJ-decomposition pieces of $E(k_1\#...\#k_n)$ are
$E(k_1),..., E(k_n)$ and $P_n\times S^1$, moreover
$E(k_1\#...\#k_n)$ is obtained by identifying $\partial E(k_i)$ and
$c_i\times S^1$ such that the meridian $m_i$ of $E(k_i)$ is
identified with $x_i\times S^1$, where $x_i$ is a point in $c_i$,
$i=1,..., n$.
\end{Lemma}

To construct our example, we need first to orient knot exteriors
and their meridians and Seifert fibers and to take a careful look
at Lemma~\ref{Schubert}.

The orientation of each knot exterior  below is induced from the
3-sphere with fixed orientation; the torus boundary of each knot
exterior has induced orientation; on each torus boundary, the
meridian and the Seifert fiber are oriented so that their product
give the orientation of the torus.

Suppose the meridian and the Seifert fiber of $E(3,2)$ have been
oriented.

\begin{Lemma}\label{ExistMap1}  (i)  The meridian and the Seifert fiber of $E(3p,2)$
can be oriented so that there is a proper map   $$\pi_p\co E(3p, 2)
\to E(3,2)$$ of degree $p$ for any odd $p$ which sends the Seifert
fiber of $E(3p, 2)$ to the $p$ times of Seifert fiber of $E(3,2)$
 and send the meridian to the meridian.

(ii) The meridian and the Seifert fiber of $\bar E(3,2)$ can be
oriented so that there is a proper degree $-1$ map  $$\bar \pi \co
\bar E(3,2) \to E(3,2)$$ which send the meridian to the meridian
and reverses the direction of the Seifert fiber.
\end{Lemma}

\begin{proof} (i) Let $A$ be a cyclic group of order $p$ acts freely on
 along the regular  Seifert fiber on $E(3p, 2)$ which induces the
identity on the base space. One can verify directly that the
quotient $E(3p, 2)/A =E(3,2)$ for odd $p$. Moreover if we lift the
orientations of the meridian and the Seifert fiber of $E(3p,2)$ to
those of $E(3p, 2)$, then the  quotient map $\pi_p : E(3p, 2) \to
E(3,2)$ meets all the conditions.

(ii) By the definition there is a proper degree $-1$ map
$$r:  \bar E(3,2) \to E(3,2)$$ induced by the mirror reflection.
Now orient the meridian and the Seifert fiber of $\bar E(3,2)$ so
that $r$ reverses the direction of meridian  and preserves the
oriented Seifert fiber. Since the trefoil knot is strongly
invertible,  there is orientation preserving involution $\tau$
which reverses both the directions of the Seifert fiber and the
meridian on $\partial E(3,2)$. Then the composition $\bar
\pi=\tau\circ r$ meets all the conditions.
\end{proof}

In the next lemma, $P_n$'s are oriented and $\partial P_n$'s have
induced orientations. The proof of the lemma is very direct.

\begin{Lemma}\label{ExistMap2} Let $d_1,..., d_n$ be integers such that $\sum d_i= 1$.
There is a proper degree one map $h(d_1, ..., d_n)\co (P_n,
c_0,\cup_{i=1}^n c_i) \to (P_1, c_0, c_1)$ such that the
restriction $h|\co c_0\to c_0$ is of degree 1 and $h|\co c_i\to
c_1$ is of degree $d_i$.
\end{Lemma}

Now we are going to construct a degree one map
$$f\co E(T(9,2)\#
\bar T(3,2)\# \bar T(3,2))\to E(3,2)$$ which we call ``folding". To
define the map, we need to present the domain and the target as
follows:
$$f\co (P_3\times S^1) \cup_{\phi_i} \sqcup_{i=1}^3 E_i\to (P_1\times S^1) \cup_\phi  E(3,2)$$
where $E_1=E({9,2}), E_2=\bar E({3,2}),  E_3=\bar E({3,2})$, and
take a careful look at $\phi_i$ and $\phi$.

First all the  meridians and the Seifert fibers of $E_i$, $i=1,2,3$,
are oriented as in Lemma \ref{ExistMap1} and all $c_i$ are oriented
as in Lemma \ref{ExistMap2}, and $S^1$ is also oriented.

Now each $\phi_i$ exactly  sends the meridian of $E_i$ to $x_i\times
S^1$. Moreover the product structure of $P_3\times S^1$ can be
chosen so that $\phi_i$ sends the Seifert fiber of $E_i$ to
$c_i\times y$, which is possible since the Seifert fiber and the
meridian of $E_i$ meets transversally in one point. The product
structure of $P_1\times S^1$ is also chosen so that $\phi$ has
similar property.

Now our map $f$ is obtained by gluing the following proper maps:

(1) $h(3, -1, -1)\times id : P_3\times S^1\to P_1\times S^1$, where
$h(3,-1,-1)$ is defined in Lemma~\ref{ExistMap2};

(2) $\pi_3: E_1 \to E(3,2)$, where  $\pi_3$ is given by
Lemma~\ref{ExistMap1} (i);

(3)  $\bar \pi: E_i\to E(3,2)$, where  $\bar \pi$ is given by
Lemma~\ref{ExistMap1} (ii), $i=2,3$.

Clearly $f$ is a proper map of degree one.

Finally we show that the map $f$ is not a de-satellization.
Otherwise there would be an essential embedded torus $T$ such that
there is a non-trivial simple closed curve $c$ which stays in the
kernel of $f_*$. Since all $E_i$ involved are small knot
exteriors, $T\subset E(k)$ must be a vertical torus in $P_3\times
S^1$, which separates $P_3\times S^1$ into two copies of
$P_2\times S^1$.
 We may that suppose $c_1$ and $c_2$ are in the same $P_2\times S^1$.
Note that $f$  send ($S^1$, $c_1$, $c_2$) of $P_2\times S^1$ to
($S^1$, $3c_1$, $-c_1$) of $P_1\times S^1$, and $c_1$ and $S^1$ form
a basis for $\pi_1(P_1\times S^1)$, one can verify directly that
there is no non-trivial simple closed curve on $T$ which stays in
the kernel of $f|: T \to P_1\times S^1$. Since $P_1\times S^1$ is
$\pi_1$-injective in $E(3,2)$, so there is no non-trivial simple
closed curve on $T$ which stays in the kernel of $f|: T \to E(3,2)$,
and we reach a contradiction. The verification of the cases that
other $c_i$ and $c_j$ are in the same $P_2\times S^1$ is similar.
\end{Example}

\end{document}